\documentclass[12]{amsart}
\usepackage{amssymb,amsfonts,amsmath}
\input{AGA.sty}
\newtheorem{lem}{{\sc Lemma}}

\begin{document}
\begin{AGApaper}

\AGAtitle{DISCRETE VERSIONS OF THE BECKMAN-QUARLES THEOREM FROM
THE DEFINABILITY RESULTS OF R. M. ROBINSON}

\AGAauthor{Apoloniusz Tyszka}

\AGAaddress{Technical Faculty,
Hugo Ko{\l}{\l}\c{a}taj University\\
Balicka 104, 30--149 Krak\'ow, Poland\\
E-mail: rttyszka@cyf-kr.edu.pl}

\begin{AGAabstract}
We present shorter proofs of the discrete versions of the
Beckman-Quarles theorem first proved in~\cite{tys}.
\end{AGAabstract}

\vskip 0.1truecm
The classical Beckman-Quarles theorem (\cite{bq}, \cite{ben})
states that if $f:{\mathbb R}^n$$\rightarrow$${\mathbb R}^n$
\mbox{$(n>1)$} preserves all unit distances then it is an isometry.
In this note we derive the discrete forms of this theorem
(\cite{tys}) from a Lemma proved there and from Robinson's
results on geometric notions which are definable in terms of the
unit distance (\cite{rob}). This new proof avoids reference to the
deeper result of H. Maehara on rigid unit-distance graphs as well
as the auxiliary results from pp.\ 128-132 in~\cite{tys}, and may
thus be considered simpler.

For the remainder of the paper we fix an integer $n>1$.
For algebraic $r>0$ let $P_r$ denote a binary predicate,
where $P_r(x,y)$ have the intended interpretation
``the distance from $x$ to $y$ is $r$''.
Let ${\mathcal L}$ be the first order language with equality in
which the only non-logical symbol is $P_1$.
We say that a $k$-ary relation $R$ in ${\mathbb R}^n$
is (existentially) definable in ${\mathcal L}$
(i.e. in terms of the unit distance using equality and logical
connectives $\neg$, $\wedge$, $\vee$, $\Rightarrow$,
$\Leftarrow$, $\Leftrightarrow$) if there exists
an (existential) ${\mathcal L}$-formula $\phi(x_1,...,x_k)$
such that
\\
\centerline{$\forall x_1,...,x_k \in {\mathbb R}^n$
$(R(x_1,...,x_k)
\Leftrightarrow \phi(x_1,...,x_k)$ holds in ${\mathbb R}^n)$.}
\vskip 0.1truecm
The following two results are special cases of
results from~\cite{rob}:
\begin{AGAtheorem}
{\em (i)} All algebraic distances in ${\mathbb R}^n$ can be defined
existentially in terms of the unit distance.\\
{\em (ii)} Let the algebraic numbers $r$ and $s$ satisfy $0<s<r$.
The local equidistance relation in ${\mathbb R}^n$:
\[s \leq d(K,L)=d(M,N) \wedge (d(K,L), d(K,M), d(K,N),
d(L,M),d(L,N),d(M,N)\leq r)\]
can be defined existentially in terms of the unit distance.
\end{AGAtheorem}

Theorem 1 in~\cite{tys} states that,
if $X,Y \in {\mathbb R}^n$ and $d(X,Y)$ is an algebraic
number then there exists a finite set
$\{X,Y\} \subseteq S_{XY} \subseteq {\mathbb R}^n$
such that\\
($\ast$) each map $f:S_{XY} \rightarrow {\mathbb R}^n$ that
preserves unit distance also preserves the distance between
$X$ and $Y$.
\par
\noindent
It means that all positive algebraic distances
in ${\mathbb R}^n$ can be defined existentially
in terms of the unit distance without using equality
and using only conjunction.

Theorem 3 in~\cite{tys} states that,
if $K,L,M,N \in {\mathbb R}^n$ and $d(K,L)=d(M,N)$ then
there exists a finite set
$\{K,L,M,N\} \subseteq C_{KLMN} \subseteq {\mathbb R}^n$,
such that\\
($\diamond$) each map $f:C_{KLMN} \rightarrow {\mathbb R}^n$
that preserves unit distance satisfies $d(f(K),f(L))=d(f(M),f(N))$.

 From item (i) of Theorem 1 follows the
existence of $S_{XY}$ with the weaker property~$(\ast)$ admitting
only injective $f:S_{XY} \rightarrow {\mathbb R}^n$ satisfying
\\
\centerline{
$\forall P,Q \in S_{XY}$ ($d(P,Q)=1
\Leftrightarrow d(f(P),f(Q))=1$).}
\\
We will denote this weaker property by $(w\ast)$.

 From item (ii) of Theorem 1 follows the existence of
$C_{KLMN}$ with the weaker property~$(\diamond)$ admitting only
injective $f:C_{KLMN} \rightarrow {\mathbb R}^n$ satisfying
\\
\centerline{
$\forall P,Q \in C_{KLMN}$ ($d(P,Q)=1
\Leftrightarrow d(f(P),f(Q))=1$).}
\\
We will denote this weaker property by $(w\diamond)$.
\\
\\
In~\cite{tys} we have proved the following:

\begin{lem}
If $X,Y \in {\mathbb R}^n$ and
$\varepsilon>0$ then there exists a finite set
$\{X,Y\} \subseteq T_{XY}(\varepsilon)\subseteq {\mathbb R}^n$,
such that each unit distance preserving mapping
$f:T_{XY}(\varepsilon) \rightarrow {\mathbb R}^n$
satisfies $|d(f(X),f(Y))-d(X,Y)| \leq\varepsilon$.
\end{lem}

If we set $\varepsilon = \frac{d(X,Y)}{2}$ in this
Lemma, we deduce that $X \neq Y$ implies $f(X) \neq f(Y)$,
if we set $\varepsilon = \frac{|d(X,Y)-1|}{2}$, we deduce that
$d(X,Y) \neq 1$ implies $d(f(X),f(Y)) \neq 1$.

Using this observation, we get the following result, which
represents a new proof of the two main results in~\cite{tys}:

\begin{AGAtheorem}
If $X,Y \in {\mathbb R}^n$ and a set
$\{X,Y\} \subseteq S_{XY} \subseteq {\mathbb R}^n$
satisfies $(w\ast)$, then
\[\widetilde{S}_{XY}:=S_{XY} \cup
\bigcup_{\stackrel{\scriptstyle P,Q \in S_{XY}}{P \neq Q}}
T_{PQ}(\frac{d(P, Q)}{2}) \cup
\bigcup_{\stackrel{\scriptstyle P,Q \in S_{XY}}{d(P,Q) \neq 1}}
T_{PQ}(\frac{|d(P,Q)-1|}{2})\]
satisfies $(\ast)$. {\em (}Since $S_{XY}$ satisfying $(w\ast)$
is known to exist by Theorem 1,
we have proved the existence of a finite set satisfying
$(\ast)${\em )}.\\
\\
If $K,L,M,N \in {\mathbb R}^n$ and a set
$\{K,L,M,N\} \subseteq C_{KLMN} \subseteq {\mathbb R}^n$ satisfies
$(w\diamond)$ then
\[\widetilde{C}_{KLMN}:=
C_{KLMN} \cup
\bigcup_{\stackrel{\scriptstyle P,Q \in C_{KLMN}}{P \neq Q}}
T_{PQ}(\frac{d(P,Q)}{2}) \cup \bigcup_{\stackrel
{\scriptstyle P,Q \in C_{KLMN}}{d(P,Q) \neq 1}}
T_{PQ}(\frac{|d(P,Q)-1|}{2})\]
satisfies $(\diamond)$. {\em (}Since $C_{KLMN}$ satisfying
$(w\diamond)$ is known to exist by Theorem 1,
we have proved the existence of a finite set satisfying
$(\diamond)${\em )}.
\end{AGAtheorem}

In~\cite{pamb} it was shown, by providing an infinitary definition
of the segment congruence relation $\equiv$ ($AB\equiv CD$ if and
only if $d(A,B)=d(C,D)$) in terms of $P_1$, that plane Euclidean
geometry over Archimedean ordered Euclidean fields (all positive
elements have square roots) can be axiomatized in the infinitary
language ${\mathcal L}_{\omega_1\omega}$.

It is worth mentioning that our results from \S 1 in~\cite{tys}
imply, for the case of Archimedean ordered Euclidean fields, more
than plain infinitary definability of $\equiv$ in terms of $P_1$.
 For such fields we have the following definition:
\[ab\equiv cd \Leftrightarrow
\bigwedge_{n=1}^\infty \bigvee_{r \in \mathbb{Q}^+}
\exists x \exists y (P_r(a,x) \wedge P_r(c,y) \wedge P_{1/n}(b,x)
\wedge P_{1/n}(d,y)),\]
in which we think of the $P_r$ and $P_{1/n}$ as abbreviations for
their positive existential definitions in terms of $P_1$,
these positive existential definitions are valid
in any Archimedean ordered Euclidean field.

\AGAreceived{02.06.2001}
\bibliographystyle{plain}
\begin{AGAbibliography}
\bibitem{bq} F.~S.~Beckman, D.~A.~Quarles, Jr., \textit{On
isometries of Euclidean spaces}, Proc. Amer. Math. Soc.~4 (1953),
810-815.
\bibitem{ben} W.~Benz,
\textit{Real geometries}, BI Wissenschaftsverlag, Mannheim, 1994.
\bibitem{pamb} V.~Pambuccian, \textit{Unit distance as single
binary predicate for plane Euclidean geometry},
Zeszyty Nauk.~Geom.~18 (1990), 5--8 (correction in:~19 (1991), 87).
\bibitem{rob} R.~M.~Robinson, \textit{Binary relations as
primitive notions in elementary geometry}, in: The Axiomatic
Method, (L.~Henkin, P.~Suppes and A.~Tarski (eds.))
North-Holland, Amsterdam, 1959, 68--85.
\bibitem{tys} A.~Tyszka, \textit{Discrete versions of the
Beckman-Quarles theorem},
Aequationes Math.~59 (2000), 124--133.
\end{AGAbibliography}
\end{AGApaper}
\end{document}